\documentclass[10pt,a4paper]{article}
\usepackage{amssymb}
\usepackage{latexsym}
\newtheorem{theorem}{Theorem}[section]

\newtheorem{lemma}{Lemma}[section]
\newtheorem{proposition}{Proposition}[section]
\newtheorem{remark}{Remark}[section]

\newtheorem{corollary}{Corollary}[section]
\def\proof{\mbox {\it \textbf{Proof.}~~}}
\usepackage{amsmath}
\numberwithin{equation}{section}
\begin{document}
\title{{\bf\Large  An autonomous Kirchhoff-type equation with general nonlinearity in $\mathbb{R}^N$}}
\author{{\small Sheng-Sen Lu}
\vspace{1mm}\\
{\it\small Center for Applied Mathematics}, {\it\small Tianjin University}\\
{\it\small Tianjin, 300072, PR China}\\
{\it\small and}\\
{\it\small Chern Institute of Mathematics and LPMC}, {\it\small Nankai University}\\
{\it\small Tianjin, 300071, PR China}\\
{\it\small e-mail: lushengsen@mail.nankai.edu.cn}}\vspace{1mm}
\date{}
\maketitle
\begin{center}
{\bf\small Abstract}
\vspace{3mm}
\hspace{.05in}\parbox{4.5in}
{{\small We consider the following autonomous Kirchhoff-type equation
\begin{equation*}
-\left(a+b\int_{\mathbb{R}^N}|\nabla{u}|^2\right)\Delta u= f(u),~~~~u\in H^1(\mathbb{R}^N),
\end{equation*}
where $a\geq0,b>0$ are constants and $N\geq1$. Under general Berestycki-Lions type assumptions on the nonlinearity $f$, we establish the existence results of a ground state and multiple radial solutions for $N\geq2$, and obtain a nontrivial solution and its uniqueness, up to a translation and up to a sign, for $N=1$. The proofs are mainly based on a rescaling argument, which is specific for the autonomous case, and a new description of the critical values in association with the level sets argument.}}
\end{center}
\noindent
{\it \small 2010 Mathematics Subject Classification}: {\small 35J15, 35J60}.\\
{\it \small Key words}: {\small Kirchhoff-type equation, Ground state, Multiplicity of solutions.}

\section{Introduction and main results}
The main concern of this paper is the following autonomous Kirchhoff-type problem with a general subcritical nonlinearity:
\begin{equation}\tag{$P$}
\left\{
\begin{aligned}
&-\left(a+b\int_{\mathbb{R}^N}|\nabla{u}|^2\right)\Delta{u}= f(u)~~\text{in}~\mathbb{R}^N,\\
&u\in H^1(\mathbb{R}^N),~~~~u\not\equiv0~~\text{in}~\mathbb{R}^N,\\
\end{aligned}
\right.
\end{equation}
where $a\geq0,b>0$ are constants, $N\geq1$, $f:\mathbb{R}\rightarrow\mathbb{R}$ is a given function that satisfies certain assumptions which will be specified later on.

In the case where $b>0$, the class of Problem $(P)$ is called of Kirchhoff type because it comes from an important application in Physic and Engineering. Indeed, if we replace $\mathbb{R}^N$ by a bounded domain $\Omega\subset\mathbb{R}^N$ in $(P)$, then we get the following Kirchhoff problem
\begin{equation*}
-\left(a+b\int_{\Omega}|\nabla{u}|^2\right)\Delta{u}= f(u),
\end{equation*}
which is related to the stationary analogue of the equation
\begin{equation*}
\rho\frac{\partial^2u}{\partial t^2}-\left(\frac{P_0}{h}+\frac{E}{2L}\int^L_0\left|\frac{\partial u}{\partial x}\right|^2\right)\frac{\partial^2u}{\partial x^2}=0
\end{equation*}
presented by G. Kirchhoff in \cite{Ki83}. Besides, $(P)$ is also called a nonlocal problem in this case because of the appearance of the term $b\left(\int_{\mathbb{R}^N}|\nabla{u}|^2\right)\Delta{u}$ which implies that $(P)$ is no longer a pointwise identity. And this phenomenon provokes some mathematical difficulties which make the studies of Problem $(P)$ particularly interesting. For more mathematical and physical background, we refer readers to papers \cite{An92,Ar96,Be40,Li84,Po75} and the references therein.

In the last ten years, by classical variational methods, there are many interesting results about the existence and nonexistence of solutions, sign-changing solutions, ground state solutions, the existence of positive solutions and positive ground states, least energy nodal solutions, multiplicity of solutions, semiclassical limit and concentrations of solutions to Kirchhoff type problems, see e.g. \cite{Az12,Az13,Az11,Deng15,Fi14,Guo15,Lu15,Xie16,Ye15} and the references therein.

On the other hand, when $a=1,b=0$, there has also been a considerable amount of research on this kind of problems during the past years. The interest comes, essentially, from two reasons: one is the fact that such problems arise naturally in various branches of Mathematical Physics, indeed the solutions of $(P)$ for the case $a=1,b=0$ can be seen as solitary waves (stationary states) in nonlinear equations of the Klein-Gordon or Schr\"{o}dinger type, and the other is the lack of compactness, a challenging obstacle to the use of the variational methods in a standard way.

In the papers \cite{Be83-3,Be83-1,Be83-2,Je03}, the authors studied the case where $a=1,b=0$, namely the following autonomous Schr\"{o}dinger problem
\begin{equation}\tag{$Q$}
\left\{
\begin{aligned}
&-\Delta{v}= f(v)&\text{in}~\mathbb{R}^N,\\
&v\in H^1(\mathbb{R}^N),~~v\not\equiv0&\text{in}~\mathbb{R}^N,
\end{aligned}
\right.
\end{equation}
under the following assumptions on the nonlinearity $f$:
\begin{itemize}
  \item[$(f_1)$] $f\in C(\mathbb{R},\mathbb{R})$ is continuous and odd.
  \item[$(f_2)$] $-\infty<\underset{t\rightarrow0}{\liminf}\frac{f(t)}{t}\leq\underset{t\rightarrow0}{\limsup}\frac{f(t)}{t}=-\omega<0$ for $N\geq3$.

      $\underset{t\rightarrow0}{\lim}\frac{f(t)}{t}=-\omega\in(-\infty,0)$ for $N=1,2$.
  \item[$(f_3)$] When $N\geq3$, $\underset{t\rightarrow\infty}{\lim}\frac{f(t)}{|t|^{\frac{N+2}{N-2}}}=0$.

  When $N=2$, for any $\alpha>0$
      \begin{equation*}
      \underset{t\rightarrow\infty}{\lim}\frac{f(t)}{e^{\alpha t^2}}=0.
      \end{equation*}
  \item[$(f_4)$] Let $F(t):=\int^t_0f(\tau)d\tau$. When $N\geq2$, there exists $\zeta>0$ such that $F(\zeta)>0$.

  When $N=1$, there exists $\zeta>0$ such that
  \begin{equation*}
  F(t)<0,~\forall~t\in (0,\zeta),~~F(\zeta)=0~~\text{and}~~f(\zeta)>0.
  \end{equation*}
\end{itemize}

Under the conditions $(f_1)-(f_3)$, it is not difficult to see that the natural functional corresponding to $(Q)$:
\begin{equation*}
S(v)=\frac{1}{2}\int_{\mathbb{R}^N}|\nabla{v}|^2-\int_{\mathbb{R}^N}F(v)
\end{equation*}
is well-defined on $H^1(\mathbb{R}^N)$ and of class $C^1$. With the aid of variational methods and critical points theory, Berestycki-Lions and Berestycki-Gallouet-Kavian established the existence results of a ground state, namely a nontrivial solution which minimizes the action functional $S$ among all the nontrivial solutions, and infinitely many bound state solutions of $(Q)$ in \cite{Be83-1,Be83-2} for $N\geq3$ and in \cite{Be83-3} for $N=2$ respectively. For the case of $N=1$, Jeanjean-Tanaka proved the existence of a nontrivial solution to $(Q)$ and the uniqueness of the nontrivial solution, up to a translation and up to a sign, in \cite{Je03}. To be more precise, they obtained the following theorem:
\begin{theorem}
Suppose that $N\geq1$ and $f$ satisfies $(f_1)-(f_4)$. Then, when $N\geq2$, Problem $(Q)$ possesses an infinite sequence of distinct solutions $\{v_k\}$ with the following properties:
\begin{description}
  \item[$~~(i)$] $v_k$ is spherically symmetric and of class $C^2$ on $\mathbb{R}^N$.
  \item[$~(ii)$] $\underset{k\rightarrow+\infty}{\lim}S(v_k)=+\infty$.
  \item[$(iii)$] $v_1$ is a positive ground state solution to $(Q)$ and decreases with respect to $r=|x|$.
\end{description}
When $N=1$, Problem $(Q)$ has a positive solution $v_1(x)$ and the set of all nontrivial solutions of $(Q)$ is $\{\pm v_1(x-c)~|~c\in \mathbb{R}\}$. In particular, all nontrivial solutions of $(Q)$ are ground state solutions, and $v_1(x)$ is spherically symmetric and decreases with respect to $r=|x|$ after a suitable translation.
\end{theorem}

\begin{remark}
\begin{description}
  \item[$~(i)$] When $N=1,2$, the assumption $(f_2)$ assumed in Theorem 1.1 can actually be generalized to a general one and the conclusions of Theorem 1.1 still hold. For detailed explanations, we refer readers to Remark 5.1 in Section 5.
  \item[$(ii)$] As a consequence of Poho\u{z}aev identity, if $v$ is any nontrivial solution of $(Q)$, then
\begin{equation*}
S(v)=\frac{1}{N}\int_{\mathbb{R}^N}|\nabla v|^2>0.
\end{equation*}
Thus, $v_1$ the ground state solution of $(Q)$ has the minimal value of $\int_{\mathbb{R}^N}|\nabla v|^2$ among all the nontrivial solutions of $(Q)$ and $\int_{\mathbb{R}^N}|\nabla v_k|^2\rightarrow+\infty$ as $k\rightarrow+\infty$.
\end{description}
\end{remark}

Under the same very general hypotheses on $f$ as above, Azzollini \cite{Az13} investigated Problem $(P)$ in the case where $a,b>0$ and $N\geq3$. In that paper,  by means of a rescaling argument based on an idea of himself in the previous paper \cite{Az12}, the author found the necessary and sufficient condition on the values of the positive parameters $a$ and $b$ in order that $(P)$ admits a nontrivial solution for the case $N\geq4$. In particular, the necessary and sufficient condition for the existence of a ground state solution was given when $N=4$. In addition, when $N=3$, Azzollini \cite{Az13} established the existence result of a ground state solution for any $a,b>0$, see also \cite{Az12} for a variational proof of this existence result. We remark here that the case $a,b>0$ and $N=1,2$ has also been considered in the more recent paper \cite{Fi14} by Figueiredo, Ikoma and J\'{u}nior and, in this case, the existence of a ground state solution was shown as partial results of \cite{Fi14}. For further details, we refer readers to \cite{Az12,Az13,Fi14}.

We also would like to mention the closely related work of Azzollini, d'Avenia and Pomponio \cite{Az11} which concerns firstly the multiplicity of solutions to Problem $(P)$. To be more precise, by variational methods, the authors established an abstract multiplicity theorem for critical points of a certain suitable perturbation of $S$. As an application, in the case where $a>0$ fixed, $b>0$ and $N\geq3$, they treated Problem $(P)$ and obtained the existence result of (at least) $k$ distinct radial solutions to $(P)$ for every $k\in \mathbb{N}$ and $b\in (0,b_k)$, where $b_k>0$ is a suitable positive constant depending on $k$.

In view of \cite{Az12,Az13,Az11,Fi14} and their results, we would like to point out a few related questions which also seem interesting:
 \begin{description}
   \item[$(Q_1)$] The existence of a ground state solution is still not clearly known when $a,b>0$ and $N\geq5$ and we wonder, in the case where $a,b>0$ and $N\geq3$, whether it is possible to prove some multiple results (perhaps included in the early paper \cite{Az11}) by using the rescaling argument introduced in \cite{Az12,Az13}.
   \item[$(Q_2)$] The degenerate case in high dimensions, that is $a=0,b>0, N\geq3$, is not considered.
   \item[$(Q_3)$] No multiplicity result is available in the case where $a,b>0$ and $N=2$, and the degenerate case in low dimensions, that is $a=0,b>0,N=1,2$, is also not concerned yet.
 \end{description}

In the present paper, we are interested in answering the questions above and improving the papers \cite{Az12,Az13,Az11,Fi14} with additional results by taking full advantage of the Theorem 1.1 in association with level sets argument and more analysis. In terms of $(f_1)-(f_3)$, we conclude that the corresponding functional $\Phi$ of $(P)$ given by
\begin{equation*}
\Phi(u)=\frac{a}{2}\int_{\mathbb{R}^N}|\nabla{u}|^2+\frac{b}{4}\left(\int_{\mathbb{R}^N}|\nabla{u}|^2\right)^2-\int_{\mathbb{R}^N}F(u)
\end{equation*}
is well-defined on $H^1(\mathbb{R}^N)$ and of class $C^1$. Now, we state our main results of this paper as follows:
\begin{theorem}
Assume that $a>0$ fixed, $b>0$, $N\geq3$ and $f$ satisfies $(f_1)-(f_4)$. Then the following statements hold.
\begin{description}
  \item[$~~(i)$] If $N=3$, $(P)$ has infinitely many distinct radial solutions $\{u_k\}$ for any $b>0$. Moreover, $u_1$ is a positive ground state solution of $(P)$, $\Phi(u_1)>0$ and $\Phi(u_k)\rightarrow+\infty$ as $k\rightarrow+\infty$.
  \item[$~(ii)$] If $N=4$, for any $k\in \mathbb{N}$, there exists a constant $b_k>0$ such that $(P)$ has at least $k$ distinct radial solutions for any $b\in (0,b_k)$ and one of them is a positive ground state solution of $(P)$ with positive energy. Moreover, we have that $b_k\rightarrow 0$ as $k\rightarrow+\infty$.
  \item[$(iii)$] If $N\geq5$, there exists a constant $b^*>0$ such that $(P)$ has a nontrivial solution if and only if $b\in (0,b^*]$, and $(P)$ has a ground state solution \textbf{if and only if} $b\in (0,b^*]$. Moreover, there is a constant $b^{**}\in (0,b^*)$ such that the global infimum of $\Phi$ on $H^1(\mathbb{R}^N)$ is negative if and only if $b\in (0,b^{**})$. As a direct consequence, $\Phi$ is nonnegative on $H^1(\mathbb{R}^N)$ if and only if $b\in [b^{**},+\infty)$. Actually, $b^*$ and $b^{**}$ have the exact expressions as follows
       \begin{equation*}
         b^*=\frac{2}{N-2}\left[\frac{N-4}{(N-2)a}\right]^{\frac{N-4}{2}}\|v_1\|^{-2}_{\mathcal{D}^{1,2}},~~~~b^{**}=\frac{4}{N}\left(\frac{N-4}{Na}\right)^{\frac{N-4}{2}}\|v_1\|^{-2}_{\mathcal{D}^{1,2}},
      \end{equation*}
      where $v_1$ is the ground state solution of $(Q)$ given by Theorem 1.1. In addition, for any $k\in \mathbb{N}$, there exists a constant $b_k>0$ such that $(P)$ has at least $k$ distinct radial solutions for any $b\in (0,b_k)$ and one of them is a positive ground state solution of $(P)$. Moreover, we have that $b_k\rightarrow 0$ as $k\rightarrow+\infty$.
\end{description}
\end{theorem}

\begin{remark}
\begin{description}
  \item[$~~(i)$] The existence of a ground state solution in the case where $a,b>0$ and $N\geq5$ is clearly known now and the multiple result is also provided here when $a,b>0$ and $N\geq3$. In addition, as we can see in Section 3, the proof of Theorem 1.2 is mainly based on the non-variational method introduced in \cite{Az12,Az13} and a new description of the critical values observed firstly by us in the present paper. Thus, Theorem 1.2 answers comprehensively the question $(Q_1)$ we raised above.
  \item[$~(ii)$] Although the multiple result given by Theorem 1.2 has already partly proved by variational methods in \cite{Az11} , inspired by \cite{Az12,Az13}, we provide a non-variational proof here which is simple and fundamental. In addition, when $N\geq5$, we can actually show the existence of not only $k$ distinct radial solutions with positive energies but also $k$ distinct radial solutions with negative energies  for any $k\in \mathbb{N}$ and sufficiently small $b>0$, see Remark 3.4 below for the detailed proof. However, in this case, the critical levels of the solutions given in \cite{Az11} are all positive. For another significant difference from \cite{Az11}, see the following Item $(iii)$.
  \item[$(iii)$] It is worth pointing out that, when $N=3$, we obtain here infinitely many distinct radial solutions for any $a,b>0$. This extends the previous multiple result of Azzollini, d'Avenia and Pomponio \cite{Az11} which only claims the existence of finitely many distinct radial solutions for suitable small $b>0$ in dimension three. Motivated partly by this significant difference, after having completed the first draft of this paper, we also managed to provide a variational proof of this extended result, see our latest work \cite{Lu16}.
  \item[$(iv)$] In the previous articles \cite{Az13,Fi14}, Azzollini and Figueiredo-Ikoma-J\'{u}nior have already shown certain results which are closely related to the existence of critical value $b^*$ given as above. While, to the best of our knowledge, the second critical value $b^{**}$ defined as above is firstly observed by us in the present paper and is totally new knowledge among recent researches in Kirchhoff-type problems.
\end{description}
\end{remark}

\begin{theorem}
Assume that $a=0$, $b>0$, $N\geq3$ and $f$ satisfies $(f_1)-(f_4)$. Then the following statements hold.
\begin{description}
  \item[$~~(i)$] If $N=3$, $(P)$ has infinitely many distinct radial solutions $\{u_k\}$ for any $b>0$. Moreover, $u_1$ is a positive ground state solution of $(P)$, $\Phi(u_1)>0$ and $\Phi(u_k)\rightarrow+\infty$ as $k\rightarrow+\infty$.
  \item[$~(ii)$] If $N=4$, $(P)$ has a nontrivial solution if and only if there exists a nontrivial solution $v$ of $(Q)$ such that $1=b\|v\|^2_{\mathcal{D}^{1,2}}$. In particular, there is no nontrivial solution for $b\in (\|v_1\|^{-2}_{\mathcal{D}^{1,2}},+\infty)$, where $v_1$ is the ground state solution of $(Q)$ given by Theorem 1.1. Suppose that $u$ is a nontrivial solution of $(P)$ for some $b\in (0,\|v_1\|^{-2}_{\mathcal{D}^{1,2}}]$, then $\Phi(u)=0$. As a direct consequence, $u$ any nontrivial solution of $(P)$ is a ground state solution of $(P)$. In addition, there exists a positive sequence $\{b_k\}$ such that, in the case where $b= b_k$, $(P)$ has \textbf{uncountably many distinct ground state solutions} $\{u_{\lambda}\}_{\lambda>0}$ which may change sign. Moreover, $u_{\lambda}$ is spherically symmetric with the properties
      \begin{equation*}
      \|u_{\lambda}\|\rightarrow +\infty~~\text{as}~\lambda\rightarrow0^+~~~~\text{and}~~~~\|u_{\lambda}\|\rightarrow 0~~\text{as}~\lambda\rightarrow+\infty,
      \end{equation*}
      and $b_k\rightarrow 0$ as $k\rightarrow+\infty$.
  \item[$(iii)$] If $N\geq5$, $\Phi$ is bounded from below and coercive with respect to $H^1(\mathbb{R}^N)$ norm, and attains the global infimum. In addition, $(P)$ has no nontrivial solutions with nonnegative energies and there exist infinitely many distinct radial solutions $\{u_k\}$ for any $b>0$. Moreover, $u_1$ is a positive ground state solution of $(P)$, $\Phi(u_k)<0$ for any $k\in \mathbb{N}$ and $\Phi(u_k)\rightarrow0$ as $k\rightarrow+\infty$. As a consequence, the global infimum of $\Phi$ on $H^1(\mathbb{R}^N)$ is negative and is achieved at $u_1$.
\end{description}
\end{theorem}

\begin{remark}
Theorem 1.3 deals with the degenerate case in high dimensions, that is $a=0,b>0, N\geq3$, and shows the existence results of a ground state solution and infinitely many distinct radial solutions. In particular, when $N=4$, we find a positive sequence $\{b_k\}^{\infty}_{k=1}$ such that $(P)$ has uncountably many ground state radial solutions for $b=b_k$. This result seems amazing, see also Corollaries 4.1 and 4.2 in Section 4. Besides, the nonexistence of nontrivial solutions is also established in dimension four when $b>0$ is large. To the best of our knowledge, all of the results above are totally new. Thus, we can say that the question $(Q_2)$ has been answered very satisfactorily.
\end{remark}

\begin{theorem}
Assume that $a\geq0,b>0$ fixed, $N=1,2$ and $f$ satisfies $(f_1)-(f_4)$.
\begin{description}
   \item[$~~(i)$] If $N=1$ then $(P)$ has a positive solution $u_1(x)$ and the set of all nontrivial solutions of $(P)$ is $\{\pm u_1(x-c)~|~c\in \mathbb{R}\}$. In particular, all nontrivial solutions of $(P)$ are ground state solutions. Moreover, after a suitable translation, $u_1(x)$ is spherically symmetric and decreases with respect to $r=|x|$.
  \item[$~(ii)$] If $N=2$ then $(P)$ has infinitely many distinct radial solutions $\{u_k\}$. Moreover, $u_1$ is a positive ground state solution of $(P)$, $\Phi(u_1)>0$ and $\Phi(u_k)\rightarrow+\infty$ as $k\rightarrow+\infty$.
\end{description}
\end{theorem}
\begin{remark}\begin{description}
                \item[$~(i)$] In order to answer the question $(Q_3)$, in Theorem 1.4, we consider the degenerate case as well as the non-degenerate case in low dimensions, and establish the multiplicity result when $N=2$ which is totally new and improves the previous paper \cite{Az11}. Motivated partly by this new knowledge, in the non-degenerate case, we gave a variational proof of this improvement in our latest work \cite{Lu16} whose preparations actually began after the completion of the first draft of the present paper.
                \item[$(ii)$] In fact, all the solutions given by Theorems 1.2-1.4 are both $C^2$ functions and they together with their derivatives up to order 2 have exponential decay at infinity for $N\geq1$.
              \end{description}
\end{remark}

The remaining part of this paper is organized as follows. In Section 2, a characterization of the solutions of $(P)$ and a new description of their energies are shown, which shall be used in an essential way in the proofs of the main results. In Sections 3-5, the proofs of Theorems 1.2-1.4 are completed respectively.

\section{Characterization and energies of the solutions}
First of all, we present the following general result which provides a characterization of the solutions of $(P)$ for us and shall play a vital role in the proofs of Theorems 1.2-1.4.
\begin{proposition}
 Assume that $a\geq0$ fixed, $b>0$ and $N\geq1$. Then $u\in H^1(\mathbb{R}^N)$ is a nontrivial solution to $(P)$ if and only if there exist $v\in H^1(\mathbb{R}^N)$ a nontrivial solution to $(Q)$ and $t>0$ such that
\begin{equation}\label{equ2.1}
t^{N-4}-a t^{N-2}=b \int_{\mathbb{R}^N}|\nabla v|^2
\end{equation}
 and $u(\cdot)=v(t\cdot)$.
\end{proposition}
\proof This proposition is a slightly modified version of Theorem 1.1 in \cite{Az13}, which relates the solutions of the nonlocal Problem $(P)$ to the solutions of the ``corresponding" Problem $(Q)$. Note that, by the elliptic regularity, every weak solution of Problem $(P)$ or $(Q)$ in $H^1(\mathbb{R}^N)$ becomes a strong solution, namely, belongs to $W^{2,p}_{\text{loc}}(\mathbb{R}^N)$ for any $1\leq p<+\infty$. Therefore, the argument of Theorem 1.1 in \cite{Az13} is also valid here and we omit the detailed proof.~~$\square$

\medskip
It is well known that solutions can be distinguished from each other by showing that their energies are different from each other. We call this procedure level sets argument. In association with Proposition 2.1, we can prove the following proposition which will be used in an essential way in the proofs of Theorems 1.2-1.4. To be more precise, the following Proposition 2.2 gives a new description of the critical values. This allows us not only to prove the multiplicity of solutions based on the level sets argument but also to be able to carry out subtle analysis in the case $N\geq5$.

\begin{proposition}
Assume that $a\geq0$ fixed, $b>0$, $N\geq1$, $u$ is a nontrivial solution of $(P)$, and $v$ and $t>0$ are related with $u$ as in Proposition 2.1. Then
\begin{equation}\label{equ2.2}
\Phi(u)=\frac{(1-at^2)\left[4-N(1-at^2)\right]}{4bNt^4}.
\end{equation}
\end{proposition}
\proof Actually, any nontrivial solution $u$ of $(P)$ satisfies the following Poho\u{z}aev identity
\begin{equation*}
\frac{N-2}{2N}a{\int_{\mathbb{R}^N}}|\nabla{u}|^2
+\frac{N-2}{2N}b\left({\int_{\mathbb{R}^N}}|\nabla{u}|^2\right)^2-\int_{\mathbb{R}^N}F(u) = 0.
\end{equation*}
As a consequence, the energy computed at any nontrivial solution $u$ of $(P)$ is
\begin{equation*}
\Phi(u)=\frac{a}{N}{\int_{\mathbb{R}^N}}|\nabla{u}|^2+b\frac{4-N}{4N}\left({\int_{\mathbb{R}^N}}|\nabla{u}|^2\right)^2.
\end{equation*}
Then, if $v$ and $t$ are related with $u$ as in Proposition 2.1, we have that
\begin{equation*}
\Phi(u)=\Phi(v(t\cdot))=\frac{a}{N}t^{2-N}{\int_{\mathbb{R}^N}}|\nabla{v}|^2+b\frac{4-N}{4N}t^{4-2N}\left({\int_{\mathbb{R}^N}}|\nabla{v}|^2\right)^2.
\end{equation*}
Equation \eqref{equ2.1} gives that
\begin{equation*}
\int_{\mathbb{R}^N}|\nabla v|^2=\frac{t^{N-4}-a t^{N-2}}{b},
\end{equation*}
and then \eqref{equ2.2} easily follows from a simple calculation.~~$\square$

\section{Proof of Theorem 1.2}

Before proving Theorem 1.2, we introduce the following lemmas which can be proved by some simple calculations.
\begin{lemma}
Assume that $a>0$ fixed, $N\geq3$ and $h(t):=t^{N-4}-a t^{N-2}$ for $t>0$.
\begin{description}
  \item[$~~(i)$] If $N=3$ then $h$ is decreasing in $t\in (0,+\infty)$ with range $(-\infty,+\infty)$.
  \item[$~(ii)$] If $N=4$ then $h$ is decreasing in $t\in (0,+\infty)$ with range $(-\infty,1)$.
  \item[$(iii)$] If $N\geq5$ then $h$ is increasing in $t\in (0,t^*)$ and decreasing in $t\in (t^*,+\infty)$ with range $(-\infty,s^*]$, where
      \begin{equation*}
      t^*:=\sqrt{\frac{N-4}{(N-2)a}}~~~~\text{and}~~~~s^*:=h(t^*)=\frac{2}{N-2}\left[\frac{N-4}{(N-2)a}\right]^{\frac{N-4}{2}}.
      \end{equation*}
      Moreover, $h>0$ for $t\in (0,a^{-\frac{1}{2}})$ and $h<0$ for $t\in(a^{-\frac{1}{2}},+\infty)$. As a consequence, for any $s\in (0,s^*)$, the equation $h(t)=s$ has exactly two positive solutions. One of them --- denoted by $\mu_1(s)$ --- belongs to $(0,t^*)$ and the other --- denoted by $\mu_2(s)$ --- belongs to $(t^*, a^{-\frac{1}{2}})$. In addition, for any $s_i>0(i=1,2)$ with $0<s_1<s_2< s^*$, we have
      \begin{equation*}
      0<\mu_1(s_1)<\mu_1(s_2)<t^*<\mu_2(s_2)<\mu_2(s_1)<a^{-{\frac{1}{2}}}.
      \end{equation*}
\end{description}
\end{lemma}
\begin{lemma}
Assume that $a,b>0$ fixed, $N\geq5$ and
\begin{equation*}
g(t):=\frac{(1-at^2)\left[4-N(1-at^2)\right]}{4bNt^4}
\end{equation*}
for $t\in(0,a^{-\frac{1}{2}})$. Then
$g$ is increasing in $t\in (0,t^*)$ and decreasing in $t\in (t^*,a^{-\frac{1}{2}})$. Moreover, $g<0$ for $t\in(0,t^{**})$ and $g>0$ for $t\in (t^{**},a^{-\frac{1}{2}})$, where
 \begin{equation*}
      t^*:=\sqrt{\frac{N-4}{(N-2)a}}~~~~\text{and}~~~~t^{**}:=\sqrt{\frac{N-4}{Na}}.
      \end{equation*}
\end{lemma}

We now employ Propositions 2.1 and 2.2 in association with Theorem 1.1 to prove Theorem 1.2 based on the level sets argument.

\bigskip
\noindent
\textbf{Proof of Theorem 1.2.}~~$(i)$ Assume that $N=3$ and $a,b>0$ fixed. In terms of Item $(i)$ of Lemma 3.1 and Proposition 2.1, for $v$ any nontrivial solution of $(Q)$ there exists a unique positive value
\begin{equation*}
t=\psi(\|v\|_{\mathcal{D}^{1,2}}):=\frac{\sqrt{4a+b^2\|v\|^4_{\mathcal{D}^{1,2}}}-b\|v\|^2_{\mathcal{D}^{1,2}}}{2a}\in (0,a^{-\frac{1}{2}})
\end{equation*}
such that $u(\cdot):=v(t\cdot)$ is a nontrivial solution of $(P)$. Moreover, in association with Proposition 2.2, we have that
\begin{equation*}
\Phi(u)=\varphi(\|v\|_{\mathcal{D}^{1,2}}):=\frac{a^2\left(2\sqrt{4a+b^2\|v\|^4_{\mathcal{D}^{1,2}}}-b\|v\|^2_{\mathcal{D}^{1,2}}\right)\|v\|^2_{\mathcal{D}^{1,2}}}{3\left(\sqrt{4a+b^2\|v\|^4_{\mathcal{D}^{1,2}}}-b\|v\|^2_{\mathcal{D}^{1,2}}\right)^2}>0.
\end{equation*}
An elementary computation shows that $\varphi$ is increasing in the value of the $\mathcal{D}^{1,2}(\mathbb{R}^N)$ norm of $v$ and $\varphi(\|v\|_{\mathcal{D}^{1,2}})\rightarrow+\infty$ as $\|v\|_{\mathcal{D}^{1,2}}\rightarrow+\infty$. In association with Theorem 1.1 and Item $(ii)$ of Remark 1.1, we conclude
\begin{equation}\label{equ3.1}
\varphi(\|v\|_{\mathcal{D}^{1,2}})\geq \varphi(\|v_1\|_{\mathcal{D}^{1,2}})>0,~~~~\underset{k\rightarrow+\infty}{\lim}\varphi(\|v_k\|_{\mathcal{D}^{1,2}})=+\infty.
\end{equation}
Without loss of generality, we may assume that $\|v_{k+1}\|_{\mathcal{D}^{1,2}}>\|v_k\|_{\mathcal{D}^{1,2}}$ for any $k\in \mathbb{N}$. Let $u_k(\cdot):=v_k(t_k\cdot)$, where $t_k:=\psi(\|v_k\|_{\mathcal{D}^{1,2}})>0$. Then, $\{u_k\}$ are the desired solutions in terms of Proposition 2.1 and \eqref{equ3.1}, and Item $(i)$ of Theorem 1.2 follows.

\bigskip
$(ii)$ Assume that $N=4$ and $v$ any nontrivial solution of $(Q)$ . In terms of Item $(ii)$ of Lemma 3.1 and Proposition 2.1, we conclude that $(P)$ has a nontrivial solution if and only if $b\in (0,\|v\|^{-2}_{\mathcal{D}^{1,2}})$, and for fixed $b\in (0,\|v\|^{-2}_{\mathcal{D}^{1,2}})$ there exists a unique positive value
\begin{equation*}
t=\sqrt{\frac{1-b\|v\|^2_{\mathcal{D}^{1,2}}}{a}}\in (0,a^{-\frac{1}{2}})
\end{equation*}
such that $u(\cdot):=v(t\cdot)$ is a nontrivial solution of $(P)$. Moreover, in association with Proposition 2.2, we have that
\begin{equation*}
\Phi(u)=\frac{a^2 \|v\|^2_{\mathcal{D}^{1,2}}}{4\left(1-b\|v\|^2_{\mathcal{D}^{1,2}}\right)}>0.
\end{equation*}
Suppose that $\bar{v}$ is another nontrivial solution of $(Q)$ with $\|\bar{v}\|_{\mathcal{D}^{1,2}}>\|v\|_{\mathcal{D}^{1,2}}$ and $b\in (0,\|\bar{v}\|^{-2}_{\mathcal{D}^{1,2}})$, then we have
\begin{equation*}
\frac{a^2 \|\bar{v}\|^2_{\mathcal{D}^{1,2}}}{4\left(1-b\|\bar{v}\|^2_{\mathcal{D}^{1,2}}\right)}>\frac{a^2 \|v\|^2_{\mathcal{D}^{1,2}}}{4\left(1-b\|v\|^2_{\mathcal{D}^{1,2}}\right)}\geq \frac{a^2 \|v_1\|^2_{\mathcal{D}^{1,2}}}{4\left(1-b\|v_1\|^2_{\mathcal{D}^{1,2}}\right)}>0
\end{equation*}
by using Item $(ii)$ of Remark 1.1 and an elementary computation. Without loss of generality, we may assume that $\|v_{j+1}\|_{\mathcal{D}^{1,2}}>\|v_j\|_{\mathcal{D}^{1,2}}$ for any $j\in \mathbb{N}$. For any given $k\in \mathbb{N}$, let
\begin{equation*}
b_k:=\|v_k\|^{-2}_{\mathcal{D}^{1,2}}>0, ~~~~u_j(\cdot):=v_j(t_j\cdot),
\end{equation*}
 where $t_j:=\sqrt{\frac{1-b\|v_j\|^2_{\mathcal{D}^{1,2}}}{a}}>0$ with fixed $b\in(0,b_k)$ and $j=1,2,\cdots,k$. It is not difficult to see that $\{u_j\}^k_{j=1}$ are $k$ distinct radial solutions of $(P)$ and $u_1$ is a positive radial ground state solution of $(P)$ with $\Phi(u_1)>0$ for any $b\in (0,b_k)$. Thus, the proof of Item $(ii)$ of Theorem 1.2 finishes.

\bigskip
$(iii)$ Assume that $N\geq5$ and $b>0$. The proof of Item $(iii)$ will be given in three steps.

\bigskip
\noindent
\textbf{Step 1.} \textit{Proof of the first conclusion in Item $(iii)$}

Proposition 2.1 and Item $(iii)$ of Lemma 3.1 imply that $(P)$ has a nontrivial solution if and only if
\begin{equation*}
b \|v\|^2_{\mathcal{D}^{1,2}}\leq s^*,
\end{equation*}
where $v$ is a nontrivial solution of $(Q)$ and $s^*>0$ is given in Item $(iii)$ of Lemma 3.1. In consideration of Item $(ii)$ of Remark 1.1, it is easy to conclude further that $(P)$ has a nontrivial solution if and only if
\begin{equation*}
b \leq b^*:=s^*\|v_1\|^{-2}_{\mathcal{D}^{1,2}},
\end{equation*}
where $v_1$ is the ground state solution of $(Q)$ given by Theorem 1.1.

We shall prove that either $u^1(\cdot):=v_1(t^1\cdot)$ or $u^2(\cdot):=v_1(t^2\cdot)$ is a ground state solution of $(P)$ for the case $b\in(0,b^*)$ and $\bar{u}:=v_1(t^*\cdot)$ is a ground state solution of $(P)$ for the case $b=b^*$ , where $t^j:=\mu_j\left(b\|v_1\|^2_{\mathcal{D}^{1,2}}\right)>0~(j=1,2)$ and $t^*>0$ are given by Item $(iii)$ of Lemma 3.1. Then the proof of the first conclusion in Item $(iii)$ of Theorem 1.2 finishes.

Actually, for the case $b\in(0,b^*)$, $u^j(j=1,2)$ defined as above are nontrivial solutions of $(P)$ in view of Proposition 2.1, and so we just need to show that either $u^1$ or $u^2$ has the minimal energy among all the nontrivial solutions of $(P)$. Assume that $u$ is any nontrivial solution of $(P)$, and $v$ and $t>0$ are related with $u$ as in Proposition 2.1, then Item $(ii)$ of Remark 1.1 gives that $\|v\|_{\mathcal{D}^{1,2}}\geq \|v_1\|_{\mathcal{D}^{1,2}}$ which implies $t\in[t^1,t^2]$ in terms of Item $(iii)$ of Lemma 3.1, and then
\begin{equation*}
\Phi(u)=g(t)\geq\min\{g(t^1),g(t^2)\}=\min\{\Phi(u^1),\Phi(u^2)\}
\end{equation*}
by taking advantage of Proposition 2.2 and Lemma 3.2. Thus, either $u^1$ or $u^2$ has the minimal energy among all the nontrivial solutions of $(P)$. In fact, in this case, we can conclude further that not $u^2$ but $u^1$ becomes the ground state solution, see Remark 3.2 below.

Similarly, we can prove that $\bar{u}:=v_1(t^*\cdot)$ is a ground state solution of $(P)$ for the case $b=b^*$.

\bigskip
\noindent
\textbf{Step 2.} \textit{Proof of the second conclusion in Item $(iii)$}

We first prove the ``if" part. Assume that $b\in(0,b^{**})$, then
\begin{equation*}
b\|v_1\|^2_{\mathcal{D}^{1,2}}<b^{**}\|v_1\|^2_{\mathcal{D}^{1,2}}=h(t^{**}),
\end{equation*}
which implies that $\underline{t}:=\mu_1(b\|v_1\|^2_{\mathcal{D}^{1,2}})\in(0,t^{**})$ by applying Item $(iii)$ of Lemma 3.1. Taking advantage of Proposition 2.2, the fact $\underline{t}\in(0,t^{**})$ and Lemma 3.2, we have
\begin{equation*}
\Phi(\underline{u})=g(\underline{t})<g(t^{**})=0,
\end{equation*}
where $\underline{u}(\cdot):=v_1(\underline{t}\cdot)$. Thus, the global infimum of $\Phi$ on $H^1(\mathbb{R}^N)$ is negative if $b\in (0,b^{**})$.

Now we prove the ``only if" part. Assume that the global infimum of $\Phi$ on $H^1(\mathbb{R}^N)$ is negative. Recall a fact in \cite{Az13} that $\Phi$ is bounded below and attains the global infimum if $N\geq5$. Then there exists a nontrivial minimizer $u$ which is also a nontrivial solution of $(P)$. Propositions 2.1 and 2.2 imply that there exist $v$ a nontrivial solution of $(Q)$ and $t>0$ such that
\begin{equation*}
h(t)=b\|v\|^2_{\mathcal{D}^{1,2}}, ~~t\in(0,a^{-\frac{1}{2}})~~\text{and}~~ g(t)=\Phi(u)<0.
\end{equation*}
From the fact that $t\in(0,a^{-\frac{1}{2}})$ and $g(t)<0$ and Lemma 3.2, we have $t<t^{**}$ which implies that
\begin{equation*}
b\|v\|^2_{\mathcal{D}^{1,2}}=h(t)<h(t^{**})=b^{**}\|v_1\|^2_{\mathcal{D}^{1,2}}
\end{equation*}
by taking advantage of Item $(iii)$ of Lemma 3.1, that is
\begin{equation*}
b<b^{**}\|v_1\|^2_{\mathcal{D}^{1,2}}\|v\|^{-2}_{\mathcal{D}^{1,2}}.
\end{equation*}
Item $(ii)$ of Remark 1.1 gives that $\|v_1\|^2_{\mathcal{D}^{1,2}}\|v\|^{-2}_{\mathcal{D}^{1,2}}\leq1$, thus $b\in(0,b^{**})$.

\bigskip
\noindent
\textbf{Step 3.} \textit{Proof of the last conclusion in Item $(iii)$.}

Without loss of generality, we may assume that $\|v_{j+1}\|_{\mathcal{D}^{1,2}}>\|v_j\|_{\mathcal{D}^{1,2}}$ for any $j\in \mathbb{N}$. For any given $k\in \mathbb{N}$, let
\begin{equation*}
b_k:=\frac{2}{N-2}\left[\frac{N-4}{(N-2)a}\right]^{\frac{N-4}{2}}\|v_k\|^{-2}_{\mathcal{D}^{1,2}}\in (0,b^*], ~~~~u_j(\cdot):=v_j(t_j\cdot),
\end{equation*}
 where $t_j:=\mu_1(b\|v_j\|^{2}_{\mathcal{D}^{1,2}})>0$ are given by Item $(iii)$ of Lemma 3.1 for fixed $b\in(0,b_k)$ and $j=1,2,\cdots,k$. It is easy to see that $\{u_j\}^k_{j=1}$ defined as above are nontrivial solutions of $(P)$. From $\|v_{j+1}\|_{\mathcal{D}^{1,2}}>\|v_j\|_{\mathcal{D}^{1,2}}(j=1,2,\cdots,k-1)$ and Item $(iii)$ of Lemma 3.1, we have
\begin{equation*}
0<t_1<\cdots<t_j<\cdots<t_k<t^*,
\end{equation*}
which implies
\begin{equation*}
\Phi(u_1)=g(t_1)<\cdots<\Phi(u_j)=g(t_j)<\cdots<\Phi(u_k)=g(t_k)<g(t^*)=\frac{a^2}{N(N-4)b}
\end{equation*}
by taking advantage of Lemma 3.2 and Proposition 2.2.

In view of Step 1 above and Remark 3.2 below, it is easy to conclude that $u_1$ is a positive ground state solution of $(P)$. Thus, $b_k$ is the desired positive constant, and then Step 3 finishes.

Summing up the above, the proof of Theorem 1.2 is completed.~~$\square$

\medskip
We would like to point out that, when $N\geq5$ and $a,b>0$, $(P)$ can also be regarded as a one-parameter problem with respect to the positive constant $a$. Under this viewpoint, similarly to the proof of Item $(iii)$ of Theorem 1.2 above, we can also establish the following analogue result.

\begin{theorem}
Assume that $b>0$ fixed, $a>0$, $N\geq5$ and $f$ satisfies $(f_1)-(f_4)$. Then there exists a constant $a^*>0$ such that $(P)$ has a nontrivial solution if and only if $a\in (0,a^*]$, and $(P)$ has a ground state solution if and only if $a\in (0,a^*]$. Moreover, there is a constant $a^{**}\in (0,a^*)$ such that the global infimum of $\Phi$ on $H^1(\mathbb{R}^N)$ is negative if and only if $a\in (0,a^{**})$. As a direct consequence, $\Phi$ is nonnegative on $H^1(\mathbb{R}^N)$ if and only if $a\in [a^{**},+\infty)$. Actually, $a^*$ and $a^{**}$ have the exact expressions as follows
\begin{equation*}
a^*=\frac{N-4}{N-2}\left[\frac{2}{(N-2)b\|v_1\|^2_{\mathcal{D}^{1,2}}}\right]^{\frac{2}{N-4}},~~~~a^{**}=\frac{N-4}{N}\left(\frac{4}{Nb\|v_1\|^2_{\mathcal{D}^{1,2}}}\right)^{\frac{2}{N-4}},
\end{equation*}
where $v_1$ is the ground state solution of $(Q)$ given by Theorem 1.1. In addition, for any $k\in \mathbb{N}$, there exists a constant $a_k>0$ such that $(P)$ has at least $k$ distinct radial solutions for any $a\in (0,a_k)$ and one of them is a positive ground state solution of $(P)$. Moreover, we have that $a_k\rightarrow 0$ as $k\rightarrow+\infty$.
\end{theorem}

Finally, we finish this section with some interesting and important comments.
\begin{remark}
It is easy to see that $\|u^1\|_{\mathcal{D}^{1,2}}>\|u^2\|_{\mathcal{D}^{1,2}}$, where $u^j(j=1,2)$ are solutions of $(P)$ given in the Step 1 above. Thus, when $N\geq5$, $a>0$ and $b\in (0,b^*)$, we can get at least two distinct positive radial solutions of Problem $(P)$ even for the prototypical nonlinearity $f(u)=- u+|u|^{p-2}u$ with $p\in (2,2^*)$. This fact stands in sharp contrast to the classical result of Kwong \cite{Kw91} which claims the uniqueness of positive radial solutions to Problem $(Q)$ with $f(v)=- v+|v|^{p-2}v$. We call this interesting difference between the cases $b>0$ and $b=0$ \textbf{uniqueness ``breaking" phenomenon}, which is induced by the nonlocal term and is new knowledge among recent researches in Kirchhoff-type problems. Here, we also remind the reader to note that the existence of two distinct positive radial solutions has actually been shown in \cite{Az13} for suitable $a,b>0$ when $N\geq5$, but there is no any further comments on this interesting result.
\end{remark}

\begin{remark}
Actually, for the case $b\in (0,b^*)$ in the Step 1 above, we can conclude further that $\Phi(u^1)<\Phi(u^2)$ which implies that $u^1$ is the desired ground state solution. For reader's convenience, we give the detailed proof here. The proof will split into two cases according as $b\in (0,b^{**}]$ or $b\in(b^{**},b^*)$.

\bigskip
\noindent
\textbf{Case 1.}~When $b\in (0,b^{**}]$, there holds
$$b\|v_1\|^2_{\mathcal{D}^{1,2}}\leq b^{**}\|v_1\|^2_{\mathcal{D}^{1,2}}=h(t^{**}).$$
By Item $(iii)$ of Lemma 3.1 and the inequality above, we have that
$$t_1:=\mu_1(b\|v_1\|^{2}_{\mathcal{D}^{1,2}})\in (0,t^{**}]~~~~\text{and}~~~~t_2:=\mu_2(b\|v_1\|^{2}_{\mathcal{D}^{1,2}})\in (t^*,a^{-\frac{1}{2}}).$$
Now, in terms of Proposition 2.2, we come to the conclusion that
$$\Phi(u^1)=g(t^1)\leq g(t^{**})=0<g(t^2)=\Phi(u^2).$$

\smallskip
\noindent
\textbf{Case 2.}
~In the remaining case $b\in(b^{**},b^*)$, the proof is different and technical. Before going further, we need to establish an auxiliary result first.

Let $c^*:=g(t^*)>0$. Then, for any $c\in (0,c^*)$, the equation $g(t)=c$ has exactly two positive solutions:
$$\tau_c:=\sqrt{\frac{(N-2)a-2\sqrt{a^2-N(N-4)bc}}{N(a^2+4bc)}}\in (t^{**},t^*)$$
and
$$\tau^c:=\sqrt{\frac{(N-2)a+2\sqrt{a^2-N(N-4)bc}}{N(a^2+4bc)}}\in (t^*, a^{-\frac{1}{2}}).$$
Let
\begin{equation*}
\gamma(c):=
\left\{
\begin{aligned}
 &h(\tau^c)-h(\tau_c),&~~~~& \text{for}~c\in (0,c^*),\\
 &0,&&\text{for}~c=c^*.
\end{aligned}
\right.
\end{equation*}
It is not difficult to see that $\gamma\in C(0,c^*]\cap C^1(0,c^*)$ and
$$\gamma'(c)=Nb\left[(\tau^c)^N-(\tau_c)^N\right]>0,~~~~\forall~ c\in (0,c^*).$$
Thus, we have that $\gamma(c)<\gamma(c^*)=0$, i.e.
\begin{equation}\label{equ3.2}
h(\tau^c)<h(\tau_c),~~~~\forall~ c\in (0,c^*).
\end{equation}

 Now, based on the conclusion above, we can complete the proof when $b\in(b^{**},b^*)$. Actually, we assume by contradiction that $\Phi(u^1)\geq \Phi(u^2)$ for some positive parameter $b\in(b^{**},b^*)$, i.e.
\begin{equation}\label{equ3.3}
g(t^1)\geq g(t^2),
\end{equation}
where
\begin{equation}\label{equ3.4}
t^1:=\mu_1\left(b\|v_1\|^2_{\mathcal{D}^{1,2}}\right)\in (t^{**},t^*),~~~~ t^2:=\mu_2\left(b\|v_1\|^2_{\mathcal{D}^{1,2}}\right)\in (t^*,a^{-\frac{1}{2}})
\end{equation}
given by Item $(iii)$ of Lemma 3.1. Let $c:=g(t^1)\in (0,c^*)$, then \eqref{equ3.3} and Lemma 3.2 give that
\begin{equation}\label{equ3.5}
\tau_c=t^1~~~~\text{and}~~~~\tau^c\in (t^*,t^2]\subset (t^*,a^{\frac{1}{2}}).
\end{equation}
It is easy to see that \eqref{equ3.4}, \eqref{equ3.2}, \eqref{equ3.5} and Item $(iii)$ of Lemma 3.1 give a contradiction:
$$h(t^2)=b\|v_1\|^2_{\mathcal{D}^{1,2}}=h(t^1)=h(\tau_c)>h(\tau^c)\geq h(t^2).$$
Thus, for any $b\in (b^{**},b^*)$, $\Phi(u^1)<\Phi(u^2)$.
\end{remark}

\begin{remark}
We know that, in the case $a,b>0$ and $N\geq5$, $\Phi$ is bounded from below and coercive with respect to $H^1(\mathbb{R}^N)$ norm, and attains the global infimum (See \cite{Az13}). Let
\begin{equation*}
\mathcal{K}:=\left\{\Phi(u)~|~\Phi'(u)=0,u\not\equiv0\right\}~~~~\text{and}~~~~m_{\inf}:=\inf_{u\in H^1(\mathbb{R}^N)}\Phi(u).
\end{equation*}
In terms of Item $(iii)$ of Theorem 1.2 and its proof, and Remark 3.2, we have now some interesting observations on $\mathcal{K}$ and $m_{\inf}$ as follows which provide us a better understanding of the functional $\Phi$.
\begin{description}
  \item[$~~(i)$] The set $\mathcal{K}$ is empty if $b>b^*$ and $\mathcal{K}=\{\frac{a^2}{N(N-4)b^*}\}$ if $b=b^*$. If $b\in(b^{**},b^*)$, then $\Phi(u^1)\in\mathcal{K}\subset [\Phi(u^1),\frac{a^2}{N(N-4)b}]$ and $\mathcal{K}\cap(\Phi(u^1),\frac{a^2}{N(N-4)b})\neq\emptyset$, where $u^1$ is given in the Step 1 above and $0<\Phi(u^1)<\frac{a^2}{N(N-4)b}$. In all the cases above, the global infimum $m_{\inf}$ always equals zero and can only be achieved by the trivial solution $u\equiv0$.
  \item[$~(ii)$] If $b=b^{**}$, then $0\in\mathcal{K}\subset[0,\frac{a^2}{N(N-4)b^{**}}]$ and $\mathcal{K}\cap (0,\frac{a^2}{N(N-4)b^{**}})\neq\emptyset$. The global infimum $m_{\inf}$ equals zero which is achieved not only by the trivial solution $u\equiv0$ but also by the positive ground state solution $u^1$.
  \item[$(iii)$] If $b\in(0,b^{**})$, then $\Phi(u^1)\in\mathcal{K}\subset[\Phi(u^1), \frac{a^2}{N(N-4)b}]$ and $\mathcal{K}\cap(0,\frac{a^2}{N(N-4)b})\neq\emptyset$, where $u^1$ is given in the Step 1 above and $\Phi(u^1)<0$. The global infimum $m_{\inf}$ is negative and can be achieved by the positive ground state solution $u^1$.
\end{description}
\end{remark}

\begin{remark} When $N\geq5$, adopting the argument in Step 3 of the proof of Theorem 1.2 $(iii)$ above, we can actually conclude further the existence of $k$ distinct solutions with positive energies and $k$ distinct solutions with negative energies for any $k\in \mathbb{N}$ and sufficiently small $b>0$. For reader's convenience, we give the detailed proof here.

Without loss of generality, we may assume that $\|v_{j+1}\|_{\mathcal{D}^{1,2}}>\|v_j\|_{\mathcal{D}^{1,2}}$ for any $j\in \mathbb{N}$. For any given $k\in \mathbb{N}$, let
\begin{equation*}
\widetilde{b}_k:=\frac{4}{N}\left[\frac{N-4}{Na}\right]^{\frac{N-4}{2}}\|v_k\|^{-2}_{\mathcal{D}^{1,2}}\in (0,b_k]\subset(0,b^{**}],
\end{equation*}
\begin{equation*}
\underline{u}_j(\cdot):=v_j(t_j\cdot)~~~~\text{and}~~~~\overline{u}_j(\cdot):=v_j(s_j\cdot),
\end{equation*}
 where $b_k>0$ is defined in Step 3 above, $t_j:=\mu_1(b\|v_j\|^{2}_{\mathcal{D}^{1,2}})$ and $s_j:=\mu_2\left(b\|v_j\|^2_{\mathcal{D}^{1,2}}\right)$ are given by Item $(iii)$ of Lemma 3.1 for fixed $b\in(0,\widetilde{b}_k)$ and $j=1,2,\cdots,k$.

 It is not difficult to conclude that $\{\underline{u}_j,\overline{u}_j\}^k_{j=1}$ defined as above are nontrivial solutions of $(P)$ and the following inequality holds
\begin{equation*}
0<t_1<\cdots<t_j<\cdots<t_k<t^{**}<t^*<s_k<\cdots<s_j<\cdots<s_1<\sqrt{a}.
\end{equation*}
And then, by taking advantage of Lemma 3.2 and Proposition 2.2, we have
\begin{equation*}
\Phi(\underline{u}_1)=g(t_1)<\cdots<\Phi(\underline{u}_j)=g(t_j)<\cdots<\Phi(\underline{u}_k)=g(t_k)<g(t^{**})=0
\end{equation*}
and
\begin{equation*}
0=g(t^{**})<\Phi(\overline{u}_1)=g(s_1)<\cdots<\Phi(\overline{u}_j)=g(s_j)<\cdots<\Phi(\overline{u}_k)=g(s_k).
\end{equation*}
\end{remark}
Thus, we complete the proof of the conclusion we claimed at the beginning of this remark.

\section{Proof of Theorem 1.3}
In this section, we shall consider the degenerate case in high dimensions, that is $a=0,b>0$ and $N\geq3$, and give a complete proof of Theorem 1.3 by taking advantage of Theorem 1.1, Item $(ii)$ of Remark 1.1 and Propositions 2.1 and 2.2. We remark here that, in the case $a=0,b>0$, \eqref{equ2.1} and \eqref{equ2.2} are reduced to
\begin{equation}\label{equ4.1}
t^{N-4}=b \|v\|^2_{\mathcal{D}^{1,2}}
\end{equation}
and
\begin{equation}\label{equ4.2}
\Phi(u)=\frac{4-N}{4bNt^4}
\end{equation}
respectively.

\bigskip
\noindent
\textbf{Proof of Theorem 1.3.}~$(i)$ The proof is similar to that of Item $(i)$ of Theorem 1.2 in Section 3 and we omit it here.

$(ii)$ The first conclusion here is a special situation of Proposition 2.1 in the case where $a=0$, $b>0$ and $N=4$. For the proof of the nonexistence result, we assume by contradiction that $(P)$ has a nontrivial solution $u$ for some positive parameter $b\in(\|v_1\|^{-2}_{\mathcal{D}^{1,2}},+\infty)$. Then the first conclusion in Item $(ii)$ of Theorem 1.3 shows that there exists a nontrivial solution $v$ of $(Q)$ such that
\begin{equation*}
1=b \|v\|^2_{\mathcal{D}^{1,2}},
\end{equation*}
which is actually impossible due to the fact that Item $(ii)$ of Remark 1.1 and $b>\|v_1\|^{-2}_{\mathcal{D}^{1,2}}$ imply
\begin{equation*}
b \|v\|^2_{\mathcal{D}^{1,2}}\geq b \|v_1\|^2_{\mathcal{D}^{1,2}}>1.
\end{equation*}

The second conclusion follows easily from Proposition 2.2, \eqref{equ4.2} and the fact $N=4$, and then the direct consequence holds.

At last, we give the proof of the last conclusion in Item $(ii)$ of Theorem 1.3. Let $b_k:=\|v_k\|^{-2}_{\mathcal{D}^{1,2}}>0$ and
\begin{equation*}
u_\lambda(\cdot):=v_k(\lambda\cdot),~~~~~~\lambda\in(0,+\infty).
\end{equation*}
It is not difficult to conclude that $\{b_k\}$ and $\{u_\lambda\}$ are the desired sequence and solutions respectively by using Proposition 2.1, the direct consequence of the second conclusion, Item $(ii)$ of Remark 1.1 and some simple calculations. Thus Item $(ii)$ of Theorem 1.3 follows.

$(iii)$ Assume that $N\geq5$, $a=0$ and $b>0$ fixed. The proof of the first conclusion is similar to that of Corollary 1.4 in \cite{Az13} and we omit it here. Next we shall prove the second conclusion, and then the consequence follows directly.

In terms of \eqref{equ4.1} and Proposition 2.1, for $v$ any nontrivial solution of $(Q)$ there exists a unique positive value
\begin{equation*}
t=b^{\frac{1}{N-4}}\|v\|^{\frac{2}{N-4}}_{\mathcal{D}^{1,2}}>0
\end{equation*}
such that $u(\cdot):=v(t\cdot)$ is a nontrivial solution of $(P)$. Moreover, in association with \eqref{equ4.2} we have that
\begin{equation*}
\Phi(u)=\beta(\|v\|_{\mathcal{D}^{1,2}}):=-\frac{N-4}{4N}b^{-\frac{4}{N-4}}\|v\|^{-\frac{8}{N-4}}_{\mathcal{D}^{1,2}}.
\end{equation*}
It is easy to see that $\beta$ is increasing in the value of the $\mathcal{D}^{1,2}(\mathbb{R}^N)$ norm of $v$, $\beta(\|v\|_{\mathcal{D}^{1,2}})<0$ and $\beta(\|v\|_{\mathcal{D}^{1,2}})\rightarrow0$ as $\|v\|_{\mathcal{D}^{1,2}}\rightarrow+\infty$. In association with Theorem 1.1 and Item $(ii)$ of Remark 1.1, we conclude
\begin{equation}\label{equ4.3}
\beta(\|v\|_{\mathcal{D}^{1,2}})\geq \beta(\|v_1\|_{\mathcal{D}^{1,2}}),~~~~\underset{k\rightarrow+\infty}{\lim}\beta(\|v_k\|_{\mathcal{D}^{1,2}})=0.
\end{equation}
Without loss of generality, we may assume that $\|v_{k+1}\|_{\mathcal{D}^{1,2}}>\|v_k\|_{\mathcal{D}^{1,2}}$ for any $k\in \mathbb{N}$. Let $u_k(\cdot):=v_k(t_k\cdot)$, where $t_k:=b^{\frac{1}{N-4}}\|v_k\|^{\frac{2}{N-4}}_{\mathcal{D}^{1,2}}>0$. Then, $\{u_k\}$ are the desired solutions in terms of Proposition 2.1, \eqref{equ4.3} and the fact that $\Phi(u_{k+1})>\Phi(u_k)(k=1,2,\cdots)$. On the other hand, the claim that $(P)$ has no nontrivial solutions with nonnegative energies follows from reduction to absurdity, and then Item $(iii)$ of Theorem 1.3 follows.~~$\square$

\medskip
To end this section, we give some special but interesting cases of the last conclusion in Item $(ii)$ of Theorem 1.3.

From the proof of the last conclusion in Item $(ii)$ of Theorem 1.3 above, we can easily get the following consequence.
\begin{corollary}
Assume that $a=0$, $b=\|v_1\|^{-2}_{\mathcal{D}^{1,2}}$, $N=4$ and $f$ satisfies $(f_1)-(f_4)$. Then Problem $(P)$ has \textbf{uncountably many distinct positive ground state solutions} $\left\{u_\lambda(\cdot):=v_1(\lambda\cdot)\right\}_{\lambda>0}$. Moreover, $u_\lambda$ is spherically symmetric with the properties
\begin{equation*}
\|u_{\lambda}\|\rightarrow +\infty~~\text{as}~\lambda\rightarrow0^+~~~~\text{and}~~~~\|u_{\lambda}\|\rightarrow 0~~\text{as}~\lambda\rightarrow+\infty.
\end{equation*}
\end{corollary}

Recall a classical result of Kwong \cite{Kw91} that signed solutions to Problem $(Q)$ with $f(v)=- v+|v|^{p-2}v$ $(2<p<2^*)$ are unique, up to a translation and up to a sign. Thus, in this case, the solutions $\{v_k\}_{k\geq2}$ given by Theorem 1.1 change sign. In terms of the proof of the last conclusion in Item $(ii)$ of Theorem 1.3 above, we can establish the following interesting result.
\begin{corollary}
Assume that $a=0$, $N=4$ and $f(u)=- u+|u|^{p-2}u$ with $p\in(2,4)$. Then, for $b=\|v_k\|^{-2}_{\mathcal{D}^{1,2}}$, Problem $(P)$ has uncountably many distinct \textbf{sign-changing ground state solutions} $\left\{u_\lambda(\cdot):=v_k(\lambda\cdot)\right\}_{\lambda>0}$. Moreover, $u_\lambda$ is spherically symmetric with the properties
      \begin{equation*}
       \|u_{\lambda}\|\rightarrow +\infty~~\text{as}~\lambda\rightarrow0^+~~~~\text{and}~~~~\|u_{\lambda}\|\rightarrow 0~~\text{as}~\lambda\rightarrow+\infty.
      \end{equation*}
\end{corollary}

\section{Proof of Theorem 1.4}

In this section, we try to deal with the low dimensions case, that is $N=1,2$. We remark that the non-degenerate case, that is $a>0$, and the degenerate case, that is $a=0$, can be treated in a uniform way here. Before proving Theorem 1.4, we introduce the following lemma which can be easily proved by some simple calculations.

\begin{lemma}
Assume that $a\geq0$ fixed, $N=1,2$ and $h(t):=t^{N-4}-a t^{N-2}$ for $t>0$.
\begin{description}
  \item[$~~(i)$] If $N=1$ and $a>0$ then $h$ is decreasing in $t\in (0,\tau)$ and increasing in $t\in (\tau,+\infty)$ with range $[s_\tau,+\infty)$, where
      \begin{equation*}
      \tau:=\sqrt{\frac{3}{a}}~~~~\text{and}~~~~s_\tau:=h(\tau)=-\frac{2\sqrt{3}}{9}a^{\frac{3}{2}}.
      \end{equation*}
      Moreover, $h>0$ for $t\in (0,a^{-\frac{1}{2}})$ and $h<0$ for $t\in(a^{-\frac{1}{2}},+\infty)$. If $N=1$ and $a=0$ then $h$ is decreasing in $t\in (0,+\infty)$ with range $(0,+\infty)$. As a consequence, the equation $h(t)=s$ has a unique positive solution --- denoted by $\mu(s)$ --- for any $s\in (0,+\infty)$ and fixed $a\geq0$.
  \item[$~(ii)$] If $N=2$ then $h$ is decreasing in $t\in (0,+\infty)$ with range $(-a,+\infty)$.
\end{description}
\end{lemma}

\bigskip
\noindent
\textbf{Proof of Theorem 1.4.}~$(i)$ Assume that $N=1$ and $a\geq0,b>0$ fixed. In terms of Theorem 1.1, Item $(i)$ of Lemma 5.1 and Proposition 2.1, there exists a unique positive value
\begin{equation*}
\nu:=\mu(b|\nabla v_1|^2_{L^2})>0
\end{equation*}
such that $u_1(\cdot):=v_1(\nu\cdot)$ is a positive nontrivial solution of $(P)$. Assume that $u$ is any nontrivial solution of $(P)$, then Proposition 2.1 gives that there exist $v$ a nontrivial solution to $(Q)$ and $t>0$ such that
$$h(t)=b|\nabla v|^2_{L^2}~~~~\text{and}~~~~u(\cdot)=v(t\cdot).$$ In view of Theorem 1.1 and Item $(i)$ of Lemma 5.1, we conclude that $v=v_1$, up to a translation and up to a sign, and $t=\nu$. Thus, $u=u_1$, up to a translation and up to a sign.

$(ii)$ The proof is similar to that of Item $(i)$ of Theorem 1.2 in Section 3 and we omit it here.~~$\square$

\begin{remark}
We would like to point out that Theorem 1.4 actually still holds if the assumption $(f_2)$ assumed in Theorem 1.4 when $N=1,2$, that is
\begin{equation}\label{equ5.1}
\underset{t\rightarrow0}{\lim}\frac{f(t)}{t}=-\omega\in(-\infty,0),
\end{equation}
is replaced by the following general one
\begin{equation}\label{equ5.2}
-\infty<\underset{t\rightarrow0}{\liminf}\frac{f(t)}{t}\leq\underset{t\rightarrow0}{\limsup}\frac{f(t)}{t}=-\omega<0.
\end{equation}
In terms of the proof of Theorem 1.4 above, we only need to show that, when $N=1,2$, the conclusions of Theorem 1.1 still hold under this generalized assumption.

In fact, in the proof of Theorem 1.1 when $N=1$ given by \cite{Je03}, \eqref{equ5.1} is just used to ensure that there exist $c,\delta>0$ such that
\begin{equation*}
tf(t)\leq-ct^2
\end{equation*}
for $t\in [-2\delta,2\delta]$. It is easy to see that this claim still holds under assumption \eqref{equ5.2}. This fact is pointed out in \cite{Fi14}.

In the proof of Theorem 1.1 for the case $N=2$ given by \cite{Be83-3}, the assumption \eqref{equ5.1} on $f$ is used in an essential way to show the Palais-Smale compactness condition for the corresponding functional under suitable constraints. It is hard to generalize \eqref{equ5.1} to the general one \eqref{equ5.2} in that argument. On the other hand, in \cite{Hi10}, Hirata, Ikoma and Tanaka gave another proof of Theorem 1.1 when $N\geq2$ by using mountain pass and symmetric mountain pass arguments to the corresponding functional $S$ directly. It is worth pointing out that the arguments explored in \cite{Hi10} are also applied to the case where $N=2$ and the assumption \eqref{equ5.1} on $f$ is replaced by the general one \eqref{equ5.2}.
\end{remark}
\begin{remark}
When $N\geq2$, under suitable additional assumptions on $f$, one can try to get nontrivial solutions of $(Q)$ with more information (for example, radial solutions with a prescribed number of nodes \cite{Mc90}, non-radial solutions \cite{Ba93,Lo04,Mu12}, etc.). This fact means that nontrivial solutions of $(P)$ with more qualitative properties can also be obtained, under some suitable conditions on the values of the nonnegative parameters $a$ and $b$ if necessary, by repeating the main arguments of this paper.
\end{remark}

\vspace{10mm}
We close with some questions closely related to the present paper that remain open.
\begin{description}
\item[~~$(i)$] It is interesting to know whether for this Berestycki-Lions type nonlinearity one can still obtain the main results of \cite{Az13} and this paper by using different new methods, e.g. variational methods, critical points theory, ODE approaches and so on. This question has already been partially answered in the non-degenerate case, see the early papers \cite{Az12,Az11}; while, to the best of our knowledge, the degenerate case of this question is totally open and seems to be the very challenging part. Here, we also remind readers to note that, after having completed the first draft of this paper, we managed to make new progress on this question in the non-degenerate case. We refer the reader to Item $(iii)$ of Remark 1.2 and Item $(i)$ of Remark 1.4 of the present paper and to our latest work \cite{Lu16} for further explanations and details respectively.

\item[~$(ii)$] Does Problem $(P)$ have infinitely many distinct solutions, up to a translation and up to a rotation, for $N\geq4$ and $a,b>0$ fixed? Theorem 1.2 of this paper gives only a partial answer, see also \cite{Az11,Lu16}.
\item[$(iii)$] The rescaling argument is used in an essential way in \cite{Az13} and this paper for the autonomous case, and does not work in the non-autonomous case, that is
    \begin{equation}\tag{$\mathcal{P}$}
    -\left(a+b\int_{\mathbb{R}^N}|\nabla{u}|^2\right)\Delta{u}= f(x,u),
    \end{equation}
    where $a\geq0$, $b>0$, $N\geq1$ and $f(x,t)$ depends on $x$. Many interesting
results for the non-autonomous Problem $(\mathcal{P})$ have been established for $N=1,2,3$, see e.g. \cite{Deng15,Guo15,Xie16,Ye15} and the references therein. By and large though, this question is still open, especially in the high dimensions case $N\geq4$.
\end{description}

\section*{Acknowledgment}

The author would like to express his sincere thanks to his advisor Professor Zhi-Qiang Wang for useful suggestions and language help. The author is also grateful to the anonymous referee for careful reading and valuable comments which led to an improvement of the first draft of this paper. This research is supported by the National Natural Science Foundation of China (No. 11271201).

\end{document}